\documentclass[11pt]{article}
\usepackage{amsmath,amsfonts,geompsfi}

\renewcommand{\bold}[1]{\medskip \noindent {\bf \boldmath #1
                        }\nopagebreak[4]}

\newcommand{\qed}{\nopagebreak[4]\hspace{.2cm} $\square$ \pagebreak[2]\medskip}
\newtheorem{theorem}{Theorem}[section]

\newcommand{\half}{{\mathbb H}}

\newcommand{\reals}{{\mathbb R}}

\newcommand{\makefig}[3]{
	\begin{figure}[htbp]
        \refstepcounter{figure}
	\label{#2}
        \begin{center}~
		#3~\\
		\medskip
                {\sf Figure \thefigure.  #1}
        \end{center}
	\medskip
	\end{figure}
}

\renewcommand{\bold}[1]{\medskip \noindent {\bf \boldmath #1
                        }\nopagebreak[4]}

\newcommand{\ital}[1]{\medskip \noindent {\em #1 }\nopagebreak[4]}

\newcommand{\bdry}{\partial}

\newcommand{\closure}{\overline}
\newcommand{\compos}{\circ}

\usepackage{amssymb}

\newcommand{\st}{\; | \;}         

\newcommand{\zbar}{{\overline{z}}}

\newcommand{\psibar}{{\overline{\psi}}}

\newcommand{\core}{\mbox{\rm core}}

\newcommand{\Hom}{\mbox{\rm Hom}}

\newcommand{\Isom}{\mbox{\rm Isom}}

\renewcommand{\mod}{\mbox{\rm mod}}
\newcommand{\Mod}{\mbox{\rm Mod}}

\newcommand{\sh}{\calS \calH}
\newcommand{\Sing}{\mbox{\rm Sing}}

\newcommand{\Teich}{\mbox{\rm Teich}}

\newcommand{\vol}{\mbox{\rm vol}}

\newtheorem{lem}[theorem]{Lemma}

\newcommand{\calC}{{\mathcal C}}

\newcommand{\calH}{{\mathcal H}}

\newcommand{\calK}{{\mathcal K}}

\newcommand{\calN}{{\mathcal N}}

\newcommand{\calS}{{\mathcal S}}

\newcommand{\calV}{{\mathcal V}}

\usepackage{euscript}

\newcommand{\eS}{{\EuScript S}}

\title{{\bf
\vspace{-.5in}
Weil-Petersson translation distance and volumes of
mapping tori
}
}

\author{Jeffrey F. Brock\thanks{Research partially by
NSF grant DMS-9971721.}}
\date{May 1, 2001}

\begin{document}
\maketitle
\begin{abstract}
Given a closed hyperbolic 3-manifold
$T_\psi$ that fibers over the circle with monodromy $\psi \colon S \to
S$, the monodromy $\psi$ determines an isometry of Teichm\"uller space
with its Weil-Petersson metric whose translation distance
$\|\psi\|_{\rm WP}$ is positive. 
We show there is a constant $K\ge 1$ depending only on the topology of $S$
so that the volume of $T_\psi$  satisfies
$\|\psi\|_{\rm WP}/K \le \vol(T_\psi) \le K \|\psi \|_{\rm WP}.$
\end{abstract}

\section{Introduction}
In this paper we generalize our study of the relation of
Weil-Petersson distance to volumes of hyperbolic 3-manifolds to the
context of those hyperbolic 3-manifolds that fiber over $S^1$.

Let $S$ be a closed surface of genus at least $2$, and let $\Mod(S)$
denote its mapping class group, or the group of isotopy classes of
orientation preserving self-homeomorphisms of $S$.
Let $\psi$ be a
{\em pseudo-Anosov} element of $\mod(S)$, i.e. an element that
preserves no isotopy class of simple closed curves on the surface
other than the trivial one.  

Elements of $\Mod(S)$ act by isometries on the Teichm\"uller space
$\Teich(S)$
with the Weil-Petersson metric.  Letting $d_{\rm WP}(.,.)$ denote
Weil-Petersson distance, let 
$$\|\psi \|_{\rm WP} = \inf_{X \in {\rm Teich}(S)} d_{\rm WP}(X, \psi(X))$$
denote the 
Weil-Petersson {\em translation distance} of $\psi$ acting as an isometry of
$\Teich(S)$.   

Choosing a particular representative $\hat{\psi} \colon S \to S$ of
$\psi$, we may form the {\em mapping torus} 
$$T_\psi = S \times [0,1]/(x,1) \sim (\hat{\psi}(x),0).$$
Thurston exhibited a hyperbolic structure on $T_\psi$ for any
pseudo-Anosov $\psi$; the  {\em hyperbolic volume} $\vol(T_\psi)$ is
then naturally associated to the mapping class $\psi$.  The question
of whether the assignment $\psi \to \vol(T_\psi)$ is a previously encountered
isotopy invariant of a homeomorphism is a natural one.

This paper relates the hyperbolic volume of $T_\psi$ to the Weil-Petersson translation
distance for $\psi$. 

\begin{theorem}
Given $S$ there is a constant $K>1$ and so that if $T_\psi$ is a
pseudo-Anosov mapping torus, we have 
$$
\frac{1}{K} \|\psi\|_{\rm WP}  \le \vol(T_\psi) 
\le K \|\psi\|_{\rm WP}.$$
\label{theorem:main}
\end{theorem}

The proof of the main theorem mirrors the proof of a similar result
for {\em quasi-Fuchsian} hyperbolic 3-manifolds \cite{Brock:wp}:
if $Q(X,Y) = \half^3/\Gamma(X,Y)$ is the hyperbolic quotient
of $\half^3$ by the Bers {\em simultaneous uniformization} of the pair
$(X,Y) \in \Teich(S) \times \Teich(S)$, then we have the following
theorem.
\begin{theorem}[Main Theorem of \cite{Brock:wp}]
Given  the surface $S$, there are constants $K_1>1$ and $K_2>0$ so that
$$\frac{d_{\rm WP}(X,Y)}{K_1} - K_2 \le
\vol(\core(Q(X,Y))) \le
K_1 d_{\rm WP}(X,Y) + K_2.$$
\label{theorem:qf}
\end{theorem}
Here, $\core(Q(X,Y))$ is the {\em convex core} of $Q(X,Y)$: the
smallest convex subset of $Q(X,Y)$ carrying its fundamental group.

The challenge of theorem~\ref{theorem:main} is primarily to relate
the techniques involved in the proof of theorem~\ref{theorem:qf}
to the setting of 3-manifolds fibering over $S^1$.  Here is a quick
sketch of the argument:
\begin{itemize}
\item The cover $M_\psi$ of $T_\psi$ corresponding to the inclusion 
$\iota \colon S \to T_\psi$ of the fiber is an infinite volume doubly
degenerate manifold homeomorphic to $S \times \reals$ and invariant by 
a translational isometry $\Psi \colon M_\psi \to M_\psi$, for which
$\Psi_*  = \psi_*$.
\item We control the volume of a fundamental domain $D$ for the action
of $\Psi$ from below by estimating up to bounded error a lower
bound for the number of bounded length closed geodesics in $D$.  Since 
$D$ is not convex, there is potential spill-over of volume into other
translates of $D$; this is rectified by a limiting argument.
\item We control the volume of $T_\psi$ from above by
building a model manifold $N_\Delta \cong T_\psi$ built out of
3-dimensional tetrahedra, and a degree-1 homotopy equivalence $f
\colon N_\Delta \to T_\psi$ that is {\em simplicial} (the lift
$\widetilde{f} \colon
\widetilde{N_\Delta} \to 
\half^3$ sends each tetrahedron to the convex hull of the image of
its vertices).  

A priori bounds on the volume of a tetrahedron in $\half^3$ give an
estimate on the total volume of the image, and by a spinning trick (as
in \cite{Brock:wp}) we may homotope $f$ to force all but 
$k \| \psi \|_{\rm WP}$ to have arbitrarily small volume.
\end{itemize}

The constructions in each case are directly those of \cite{Brock:wp},
to which we refer often in the interest of brevity.  The paper will
introduce vocabulary in section~\ref{section:vocabulary} and prove the 
main theorem in section~\ref{section:proof}.

\bold{Acknowledgements.}  I would like to thank Yair Minsky and Dick
Canary for inspirational discussions on the topic of this paper.

\section{Vocabulary}
\label{section:vocabulary}
We define our terms.  Throughout, $S$ will be a closed surface of
genus at least $2$.  

\bold{Teichm\"uller space.}  The {\em Teichm\"uller space} $\Teich(S)$ is the
set of marked hyperbolic structures on $S$, namely pairs $(f,X)$ of a
hyperbolic surface $X$ and a homeomorphism $f\colon S \to X$, up to
the equivalence relation $(f,X) \sim (g,Y)$ if there is an isometry
$\phi \colon X \to Y$ so that $\phi \compos f \simeq g$.  The
Teichm\"uller space inherits a topology from the marked bi-Lipschitz
distance: 
$$d_{\rm bL}((f,X),(g,Y)) = \inf_{\phi \simeq g \compos f^{-1}} \log(
L(\phi))$$
where $$L(\phi) = \sup_{x \in X, v \in T_x(X)} |D\phi(v)|/|v|$$
denotes the best bi-Lipschitz constant for $\phi$.

Considering $X = \half/\Gamma$ as a Riemann surface, let $Q(X)$ denote
the vector space of holomorphic quadratic differentials $\phi(z) dz^2$
on $X$.  The space $Q(X)$ is naturally the holomorphic cotangent space
to $\Teich(S)$.  The Weil-Petersson Hermitian metric on $\Teich(S)$ is
induced by the inner product on $Q(X)$
$$\langle \phi, \psi \rangle_{\rm WP} \int_X \frac{\phi
\psibar}{\rho^2} dz d\zbar$$ (where $\rho(z) |dz|$ is the Poincar\'e
metric on $X$) via the duality
$(\mu,\phi) = \int_X \mu \phi$ where $\mu$, a {\em Beltrami
differential}, is a tangent vector to $\Teich(S)$.  Here, we consider
the Riemannian part $g_{\rm WP}$ of the Weil-Petersson metric, and its
associated distance function $d_{\rm WP}(.,.).$

We will often refer to a hyperbolic Riemann surface $X$, assuming an
implicit marking $(f \colon S \to X).$

\bold{Hyperbolic 3-manifolds.}  Moving up a dimension, we denote by 
$$\calV(S) = \Hom(\pi_1(S) ,\Isom^+(\half^3)) / {\rm conjugation}$$
the {\em representation variety} for $\pi_1(S)$.  
Each equivalence class $[\rho]$ lying in the subset $AH(S) \subset
\calV(S)$ consisting of representations that are discrete and faithful
determines a complete hyperbolic 3-manifold $M =
\half^3/\rho(\pi_1(S))$.  Since hyperbolic 3-manifolds are
$K(\pi,1)$s, points in $AH(S)$ are in bijection with the set of all
marked hyperbolic 3-manifolds $(f \colon S \to M)$, where $f$ is a
homotopy equivalence, modulo the equivalence relation
$(f \colon S \to M) \sim (g \colon S \to N)$ if there is an isometry
$\phi \colon M \to N$ so that $\phi \compos f \simeq g$.

\bold{Simple closed curves.} Let $\eS$ denote the set of isotopy
classes of essential simple closed 
curves on $S$.  Given two curves $\alpha$ and $\beta$ in $\eS$ their
{\em geometric intersection number} $i(\alpha,\beta)$ is obtained by
counting the minimal number of points of intersection over all
representatives of $\alpha$ and $\beta$ on $S$.

A {\em pants decomposition} $P$ of $S$ is a maximal collection of
distinct elements of $\eS$ with pairwise disjoint representatives on
$S$.  Two pants decompositions $P$ and $P'$ are related by an {\em
elementary move} if $P'$ can be obtained from $P$ by removing a curve
$\alpha \in P$ and replacing it with a curve $\beta \not= \alpha$ so
that $i(\alpha,\beta)$ is minimal among all choices for $\beta$ (see figure~\ref{figure:move}).
Denote by $C_{\bf P}(S)$ the graph whose vertices are pants
decompositions of $S$ and whose edges join pants decompositions
differing by an elementary move.
\makefig{Elementary moves on pants
decompositions.}{figure:move}{\psfig{file=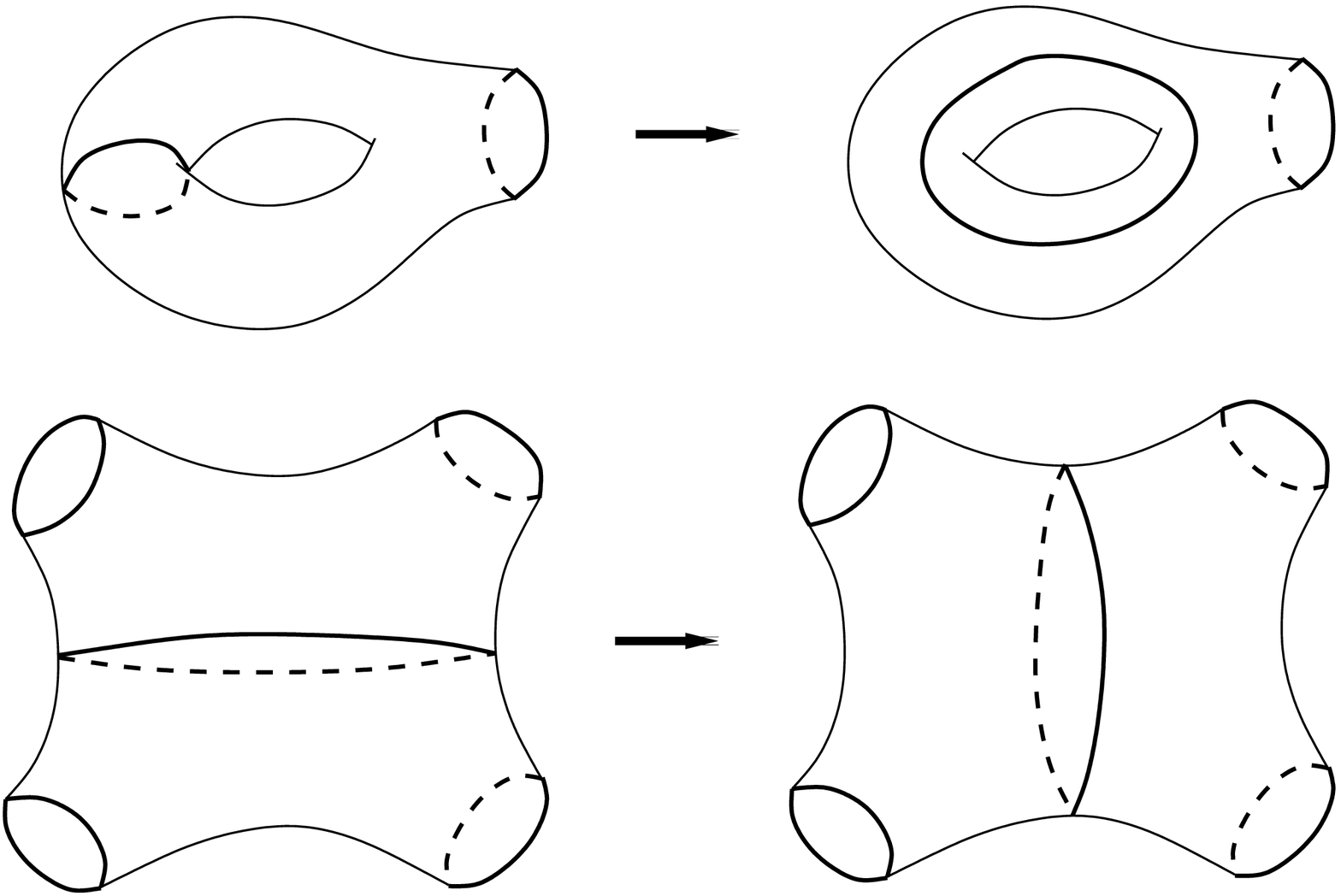,height=2in}}
Hatcher and Thurston introduced the complex $C_{\bf P}(S)$ in
\cite{Hatcher:Thurston:pants}, where they prove $C_{\bf P}(S)$ is
connected (see also \cite{Hatcher:pants}).  It therefore carries a
natural path metric, obtained by assigning each edge length 1.  We
will be interested in the distance between vertices: we let $d_{\bf
P}(P,P')$ to be the length of the minimal length path in $C_{\bf
P}(S)$ joining $P$ to $P'$.

In \cite{Brock:wp} we proved the following theorem.
\begin{theorem}
The graph $C_{\bf P}(S)$ is coarsely quasi-isometric to $\Teich(S)$
with its Weil-Petersson metric.
\label{theorem:pants:qi}
\end{theorem}
(See \cite[Thm 1.1]{Brock:wp}).
A quasi-isometry $Q \colon C_{\bf P}(S) \to \Teich(S)$ is obtained by
taking $Q$ to be any mapping for which $Q(P) = X_P$ 
lies in the set
$$V(P) = \{ X \in \Teich(S) \st \ell_X(\alpha) < L \ \ \text{for each}
\ \ \alpha \in P\}$$
where $L$ is {\em Bers' constant}, namely,
the minimal constant for which each hyperbolic surface $X \in
\Teich(S)$ admits a pants decomposition $P$ all of whose
representatives have length less than $L$ on $X$ (see
\cite{Buser:book:spectra}). 

\bold{Simplicial hyperbolic surfaces.}  Let $\Sing_k(S)$ denote the
marked {\em singular} hyperbolic structures on $S$ with at most $k$
cone singularities, namely surfaces that are hyperbolic away from at
most $k$ cone points, each with cone angle at least $2\pi$, equipped
with homeomorphisms $g \colon S \to Z$, up to marking preserving
isometry.  As with $\Teich(S)$, $\Sing_k(S)$ is topologized via 
the marked bi-Lipschitz distance, defined analogously.

Let $M$ be a complete hyperbolic 3-manifold.
A {\em simplicial hyperbolic surface} 
is a path-isometric map $h \colon Z \to M$ of a singular hyperbolic
surface into $M$ so that $h$ is totally geodesic on the faces of a
geodesic {\em triangulation} $T$ of $Z$, in the sense of
\cite{Hatcher:triangulations} (cf. \cite[Sec. 4]{Brock:wp}).  We say
$g$ is {\em adapted} to $T$.  

Given $(f \colon S \to M)$ in $AH(S)$, we denote by $\sh_k(S)$ the
{\em simplicial hyperbolic surfaces homotopic to $f$}, in other words,
those simplicial hyperbolic surfaces $h \colon Z \to M$, so that if $g
\colon S \to Z$ is the marking on $Z$, then $h \compos g$ is homotopic
to $f$.  

As in \cite{Brock:wp}, the following result of Canary (see
\cite[Sec. 5]{Canary:inj:radius}) will be
instrumental in what follows.
\begin{theorem}[Canary]
Let $M \in AH(S)$ have no accidental parabolics, and let $(g_1,Z_1)$
and $(g_2,Z_2)$ lie in $\sh_1(M)$ where $(g_1,Z_1)$ and $(g_2,Z_2)$
are adapted to $\alpha$ and $\beta$.  Then there is a continuous
family $(g_t \colon Z_t \to M) \subset \sh_2(M)$, $t \in [1,2]$.
\label{theorem:canary}
\end{theorem} 

\section{Proof of the main theorem}
\label{section:proof}
Given $\psi \in \Mod(S)$, let 
$$T_\psi = S \times I /(x,0) \sim (\psi(x),1)$$
be the mapping torus for $\psi$.  
By a theorem of Thurston (\cite{Thurston:hype2}, see also
\cite{McMullen:book:RTM}, \cite{Otal:book:fibered}), $T_\psi$ admits a
complete hyperbolic structure if and only if $\psi$ is {\em
pseudo-Anosov}: no power of $\psi$ fixes any non-peripheral isotopy
class of essential simple closed curves on $S$.

To prove theorem~\ref{theorem:main} we relate the volume of the
mapping torus $T_\psi$ to the translation distance of $\psi$ on the
pants complex $C_{\bf P}(S)$.

\bold{Proof:}{{\em (of theorem~\ref{theorem:main})}  Let $M_\psi$ be
the cover of $T_\psi$ corresponding to 
the inclusion of the fiber $\iota \colon S \to T_\psi$; $\iota$
naturally lifts to a marking $(\tilde{\iota} \colon S \to M_\psi)$ of
$M_\psi$, so $M_\psi$ is naturally an element of $AH(S)$.  

We will apply the techniques of \cite{Brock:wp} directly to larger and
larger subsets of the manifold $M_\psi$.

\ital{Bounding volume from below.}
Consider any simplicial hyperbolic surface $(h \colon Z \to M_\psi)
  \in \sh_1(M_\psi)$.  Let $(h_\psi \colon Z_\psi \to M_\psi) \in
\sh_1(M_\psi)$ be the simplicial hyperbolic surface obtained by
post-composing the image with the generating covering transformation
$\Psi \colon M_\psi \to M_\psi$ acting on $\pi_1(M_\psi)$ by $\psi_*$.

The hyperbolic surface $Z^{\rm h}$ in the same conformal class as
$Z$ differs from the hyperbolic surface $Z_\psi^{\rm h}$ by the
action of $\psi$ on $\Teich(S)$.  If $Z^{\rm h}$ lies in $V(P)$, then,
the surface $Z_\psi^{\rm h}$ lies in $V(\psi(P))$.

By an argument using
theorem~\ref{theorem:canary} in an exactly 
analogous manner to its use in \cite[Sec. 4]{Brock:wp}, 
there is a sequence 
$$g = \{P = P_0, P_1,\ldots, P_n = \psi(P)\}$$
and a continuous family $(h_t
\colon Z_t \to M_\psi) \subset \sh_k(M_\psi)$, $k =2 $, $0 \le
t \le n$, of
simplicial hyperbolic surfaces joining $Z$ and $Z_\psi$ so that
the following holds:
\begin{itemize}
\item For each $i$, the simplicial hyperbolic surface $h_i \colon Z_i
\to M_\psi$ satisfies
$\ell_{Z_i^{\rm h}}(\alpha) < L$ for each $\alpha \in P_i$.
\item Successive pants decompositions have distance $d_{\bf
P}(P_i,P_{i+1}) < D$ where $D$ depends only on $S$.
\end{itemize}

Let $\eS_g \subset \eS$ denote the set of curves 
$$\eS_g = \{ \alpha \in \eS \st \alpha \in P_i,  \ P_i \in g \}.$$
Then we phrase the following lemma of \cite{Brock:wp} (lemma 4.2) to
suit our
context:
\begin{lem}
There is a constant $K$ depending
only on $S$ so that 
$$d_{\bf P}(P,\psi(P)) \le K |\eS_g|$$ 
where $|\eS_g|$ denotes the number of elements in $\eS_g$.
\label{lemma:count}
\end{lem}

Consider the image of the homotopy $H \colon [0,n] \times S \to
M_\psi$ for which $H(t,.) = h_t \compos g_t \colon S \to M_\psi$
and $g_t \colon S \to Z_t$ is the implicit marking on $Z_t$.

Applying the argument of lemma 4.1 of \cite{Brock:wp}, there are constants
$c_1>1$ and $c_2 >0$  depending only on $S$ so that the $\epsilon$
neighborhood of the image $\calN_\epsilon( H([0,n] \times S))$ has
volume
$$\vol (\calN_\epsilon( H([0,n] \times S))) > \frac{|\eS_g|}{c_1} - c_2 \ge
\frac{d_{\bf P} (P,\psi(P))}{c_1 K} - c_2.$$

The image $H([0,n] \times S )$ is a compact subset of $M_\psi$; choose
an embedding $h \colon S \to M_\psi$ homotopic to $\iota$ so that
$h(S)$ does not intersect the $\epsilon$-neighborhood
$\calN_\epsilon(H([0,n] \times S))$ of the image $H([0,n] \times S)$.  
Let $\Psi\colon M_\psi \to M_\psi$ be the isometric covering
transformation so that 
$$\iota \compos \psi  \simeq \Psi\compos \iota.$$
Let $n_0$ be an integer so that $\Psi^{n_0} \compos h (S)$
is also disjoint from  $\calN_\epsilon(H([0,n]\times S))$ and so that
$\calN_\epsilon(H([0,n] \times S))$ lies in the compact subset $\calK$
of $M_\psi$ bounded by 
$h(S)$ and $\Psi^{n_0}\compos h(S)$.  
Note that 
$$\vol(\calK_{n_0}) = n_0 \vol(T_\psi)$$
since the region bounded by each $\Psi^i \compos h(S)$ and
$\Psi^{i+1} \compos h(S)$ is a fundamental domain for the action on
$\Psi$ on $M_\psi$.
\makefig{Bounding volume from below.}{figure:translate}{\psfig{file=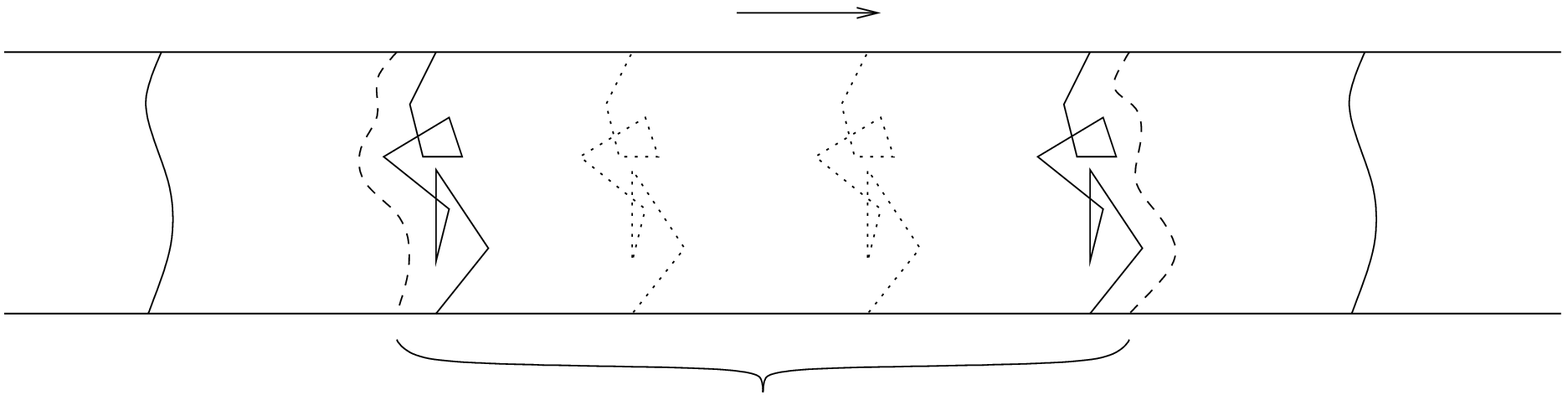,height=1.4in}}

It follows that if 
$$\calC_j = \bigcup_{i=0}^j \;\Psi^i \left( \calN_\epsilon(H([0,n] \times S))
\right)$$
then we have
$$(n_0 + j) \vol (T_\psi) \ge \vol (\calC_j) 
\ge
\frac{ j d_{\bf P} (P,\psi(P))}{c_1 K} - j c_2.$$
Taking the limit as $j$ tends to infinity, we have
$$\vol(T_\psi) \ge \frac{d_{\bf P}(P,\psi(P))}{c_1 K} - c_2.$$

Applying theorem~\ref{theorem:pants:qi},
there are constants $d_1 >1$ and $d_2 >0$ so that
$$\vol(T_\psi) \ge \frac{\|\psi\|_{\rm WP}}{d_1} - d_2.$$

But given $S$, there is an $\epsilon_S$ so that 
$$\min \{ \vol(T_\psi) \st \psi \in \Mod(S) \ \ \text{is
pseudo-Anosov} \}$$ is greater than $\epsilon_S$, 
 so there is a thus a $K_0 >1$ so that 
$$\vol(T_\psi) \ge \frac{\|\psi\|_{\rm WP}}{K_0},$$
proving one direction of the theorem.

\ital{Bounding volume from above.}
To bound $\vol(T_\psi)$ from above, 
we adapt our construction of a simplicial model 3-manifold for a
quasi-Fuchsian 3-manifold in
\cite[Sec. 5]{Brock:wp} to build a triangulated 3-manifold $N \cong
T_\psi$, and a degree one homotopy equivalence $f \colon N \to T_\psi$
that is {\em simplicial}: the lift $\widetilde{f} \colon \widetilde{N}
\to \half^3$ maps each simplex $\Delta \subset N$ to the convex hull
of the images of its vertices.

Consider the surface $Z^{\rm h}$ in $V(P)$, and let $T$ 
be a standard triangulation suited to $P$, in the sense of
\cite[Defn. 5.3]{Brock:wp}.  Let
$(g_0 \colon Z_0 \to M_\psi) \in \sh_k(M_\psi)$ be a simplicial
hyperbolic surface adapted to $T$ that realizes each $\alpha \in P$.
As in \cite{Brock:wp} we may choose $(g_0 \colon Z_0 \to M_\psi)$ so
that the vertices of $T$ map to pairs of antipodal vertices on each
geodesic $\alpha^*$; in other words the two vertices that lie in
$\alpha^*$ separate $\alpha^*$ into two segments of the same length.

Let $(g_1 \colon Z_1 \to M_\psi)$ be the simplicial hyperbolic surface
$$g_1 = \Psi\circ g_0 \circ \psi^{-1}.$$
Since $g_1 \circ \psi$ is homotopic to $\Psi\circ g_0$, the
simplicial hyperbolic surface $(g_1 \colon Z_1 \to M_\psi)$ lies in
$\sh_k(M_\psi)$.  Moreover, $g_1$ realizes the pants decomposition
$\psi(P)$ and is adapted to $\psi(T)$.

As shown in \cite{Brock:wp}, there is a triangulated model 3-manifold $N
\cong S \times I$ built out of 3-dimensional tetrahedra, together
with a simplicial mapping $g \colon N \to M_\psi$ 
with the following properties:
\begin{enumerate}
\item There is a constant $k$ so that all but $k d_{\bf P}(P,
\psi(P))$ of the tetrahedra in $N$ have the property that one edge
maps by $g$ to a geodesic $\alpha^*$, where $\alpha \in \eS_g$. 
\item If $\bdry^- N = S \times \{0\}$ and $\bdry^+ N = S \times \{1\}$, then
$g \vert_{\bdry^-N}$ factors through the simplicial hyperbolic
surface $(g_0 \colon Z_0 \to M_\psi)$, and $g \vert_{\bdry^+ N}$
factors through the simplicial hyperbolic surface $(g_1 \colon Z_1 \to
M_\psi)$.
\end{enumerate}

After adjusting $g \vert_{\bdry^+ N}$ by precomposition with an
isotopy,  we may glue $\bdry^- N $ to $\bdry^+ N$ by a homeomorphism
$h \simeq \psi$ to obtain a manifold $N_\psi \cong T_\psi$, and 
a homotopy equivalence $f_\psi \colon N_\psi \to T_\psi$.  By a
theorem of Stallings \cite{Stallings:fibers}, $f_\psi$ is homotopic to
a homeomorphism, and is thus surjective since $T_\psi$ is closed.  

There is a constant $\calV_3$ so that the volume of each tetrahedron
in $\half^3$ is less than $\calV_3$, so the volume of the image of $g$
is less than $\calV_3$ times the number of tetrahedra in $N$.  Arguing
as in \cite[Sec. 5]{Brock:wp}, we may {\em spin} the mapping $f_\psi$
by pulling the vertices $p_\alpha$ and $\closure{p}_\alpha$ around the
geodesic $\alpha^*$, keeping the map simplicial.  In the process, the
volume of the image of each tetrahedron $\Delta$ in $N$ for which $g$
sends an edge of $\Delta$ to $\alpha^*$ for some $\alpha \in \eS_g$,
can be made arbitrarily small by spinning the mapping sufficiently far
(see \cite[Lem. 5.10]{Brock:wp}).

We can, then, for each $\epsilon$ spin the mapping $f_\psi$ to a
mapping $f_\psi^\theta$ so that the image of $f_\psi^\theta$ has
volume less than  
$$k\calV_3  d_{\bf P}(P,\psi(P)) + \epsilon.$$
It follows from theorem~\ref{theorem:pants:qi} that there are
constants $e_1 >1$ and $e_2 >0$ so that 
$$\vol(T_\psi) \le e_1 d_{\rm WP}(Z, \psi(Z)) + e_2.$$
But there is a constant $C(Z)$ depending on the surface $Z$, so that
$$d_{\rm WP}(Z, \psi^n(Z)) \le C(Z) + \|\psi^n \|_{\rm WP} = C(Z) + n
\|\psi\|_{\rm WP}.$$  (The constant $C(Z)$ measures twice the distance
to the min-set for the action of the isometry $\psi$ of the
Weil-Petersson completion $\closure{\Teich(S)}$).  

Since $\vol(T_{\psi^n}) = n \vol(T_\psi)$, we have
$$n \vol(T_\psi) \le e_1 (C(Z) + n \|\psi\|_{\rm WP}) + e_2.$$
Taking limits of both sides as $n \to \infty$, we have
$$\vol(T_\psi) \le e_1 \|\psi\|_{\rm WP}.$$
Setting $K = \max\{K_0,e_1\}$ proves the theorem. \qed

\bibliographystyle{math}
\bibliography{math}

{\noindent \scriptsize \sc Math Department, University of
Chicago, 5734 S. University Ave., Chicago, IL 60637\\
email: brock@math.uchicago.edu}

\end{document}